\documentclass[12pt]{article}
\usepackage{mathrsfs}
\usepackage{stmaryrd}
\usepackage{amsfonts,amsmath,amssymb,amscd}
\usepackage{shadow}
\usepackage{graphicx}
\usepackage{color}
\usepackage{longtable}
\usepackage{pstricks,multido}
\usepackage{hyperref}
\allowdisplaybreaks

\newtheorem{thm}{Theorem}[section]
\newtheorem{lem}[thm]{Lemma}

\newtheorem{cor}[thm]{Corollary}

\newtheorem{conj}{Conjecture}[section]

\def\qed{\hfill \rule{4pt}{7pt}}

\def\pf{\noindent {\it{Proof.} \hskip 2pt}}

\numberwithin{equation}{section}

\begin{document}
\begin{center}
{\large\bf  Ordered Partitions Avoiding
a Permutation  of Length 3}
\end{center}

\begin{center}
William Y.C. Chen$^1$, Alvin Y.L. Dai$^2$, Robin D.P. Zhou$^3$

$^{1,2,3}$Center for Combinatorics, LPMC-TJKLC\\
Nankai University\\
Tianjin 300071, P.R. China

$^1$Center for Applied Mathematics\\
Tianjin University\\
Tianjin 300072, P.R. China

$^1$chen@nankai.edu.cn,
$^2$alvin@cfc.nankai.edu.cn,
$^3$robin@cfc.nankai.edu.cn

\end{center}

\begin{abstract}
 An ordered  partition of $[n]=\{1, 2, \ldots, n\}$ is
a partition whose blocks are endowed with a linear order.
Let $\mathcal{OP}_{n,k}$ be set of ordered partitions of $[n]$ with $k$ blocks
and $\mathcal{OP}_{n,k}(\sigma)$ be set of ordered partitions in $\mathcal{OP}_{n,k}$ that avoid a pattern $\sigma$. Recently, Godbole, Goyt, Herdan and Pudwell obtained formulas for the number of  ordered   partitions of $[n]$
with $3$ blocks and the number of ordered  partitions of $[n]$ with $n-1$ blocks avoiding a permutation pattern of length $3$. They showed that
$|\mathcal{OP}_{n,k}(\sigma)|=|\mathcal{OP}_{n,k}(123)|$ for any permutation $\sigma$ of length $3$, and raised the question concerning the
enumeration of  $\mathcal{OP}_{n,k}(123)$.
They also conjectured that the number of  ordered  partitions of $[2n]$
with blocks of size $2$ avoiding a permutation pattern of length $3$ satisfied a second order linear recurrence relation.
In answer to the question of Godbole, et al.,
we obtain the generating function for  $|\mathcal{OP}_{n,k}(123)|$ and we prove
 the conjecture on the recurrence relation.
\end{abstract}

\noindent {\bf Keywords}:  pattern avoidance, ordered  partition

\noindent {\bf AMS  Subject Classifications}: 05A15, 05A18

\section{Introduction}

The notion of pattern avoiding permutations was introduced by Knuth \cite{Knuth}, and it has been extensively studied. Klazar \cite{M.Klazar} initiated the study of
pattern avoiding set partitions. Further studies of pattern avoiding set partitions can be found in \cite{A.Goyt,Mansour,M.Klazar2,M.Klazar3,B.Sagan}.
Recently,  Godbole,  Goyt,  Herdan and  Pudwell \cite{Pudwell} considered pattern avoiding ordered set partitions. From now on, a partition means
a set partition. Let $[n]=\{1, 2,\ldots, n\}$.
  For a permutation pattern $\sigma$ of length $3$,
they obtained formulas for the
number of $\sigma$-avoiding  ordered   partitions of $[n]$ with $3$ blocks  and
the number of $\sigma$-avoiding ordered partitions of $[n]$ with $n-1$ blocks. They
raised the question of finding the number  of ordered partitions of $[n]$ with
$k$ blocks avoiding a permutation pattern of length $3$.

In answer to the above question, we  derive a bivariate generating function for the number of ordered  partitions of $[n]$ with
$k$ blocks avoiding a permutation pattern of length $3$. Meanwhile, we confirm
 the conjecture also posed by Godbole, et al. on the recurrence relation concerning the number of ordered partitions of $[2n]$ with blocks of size $2$ avoiding a permutation pattern of length $3$.

Let us give an overview of notation and terminology.
Let $S_n$ be the set of permutations of $[n]$.
Given a permutation $\pi=\pi_1\pi_2\cdots\pi_n$ in $S_n$ and  a permutation $\sigma=\sigma_1 \sigma_2 \cdots \sigma_k \in S_k$, where $1\leq k\leq n$, we say
 that $\pi$ contains a pattern $\sigma$  if there exists a subsequence
 $\pi_{i_1}\pi_{i_2} \cdots \pi_{i_k}$ $(1 \leq i_1 < i_2 < \cdots <i_k \leq n)$ of $\pi$ that is order-isomorphic to $\sigma$,
 in other words,  for all $l,m \in [k]$, we have
 $\pi_{i_l}<\pi_{i_m}$ if and only if $\sigma_l <\sigma_m$. Otherwise, we say that $\pi$  avoids a pattern $\sigma$, or
 $\pi$ is $\sigma$-avoiding. Let $S_n(\sigma)$ denote the set of permutations of $S_n$ that avoid a pattern $\sigma$.
 For example, $41532$ is 123-avoiding, while it
 contains a pattern $312$ corresponding to the subsequence $412$.

 A partition $\pi$ of the set $[n]$, written
 $\pi\vdash [n]$, is
a family of nonempty, pairwise disjoint subsets $B_1,B_2,\ldots, B_k$ of
$[n]$ such that $\cup_{i=1}^{k}B_i=[n]$, where each $B_{i}~(1\leq i\leq k)$ is called  a block of partition $\pi$. We write $\pi=B_1/B_2/\cdots /B_k$ and define the length of $\pi$, denoted $b(\pi)$, to be the number of blocks.
 For convenience, we may write elements of a block in
increasing order and list the blocks in the increasing order of their minimum elements.  We represent such a partition by a canonical sequence
$\rho=\rho_1\rho_2\cdots \rho_n$ with $\rho_i=j$ if $i\in B_j$.
Let $\sigma$ be a permutation of $[m]$ with $m\leq n$. We say that a partition $\pi$ contains a pattern $\sigma$ if the canonical sequence of $\pi$ contains a subsequence that is order-isomorphic to $\sigma$. Otherwise, we say $\pi$ avoids a pattern $\sigma$.
For example, given a partition $\pi=156/24/37$, we have $\pi\vdash [7]$, $b(\pi)=3$. The canonical sequence of $\pi$ is $1232113$, which is $312$-avoiding.

An ordered partition of $[n]$ is
a partition of $[n]$ whose blocks are endowed with a linear order. Let $\mathcal{OP}_{n,k}$ denote the set of  ordered partitions of $[n]$ with $k$ blocks, let  $\mathcal{OP}_n$  denote the set of ordered partitions of $[n]$, and
let $\mathcal{OP}_{[b_1,b_2,\ldots,b_k]}$ denote the set of ordered partitions   of $[b_1+b_2+\cdots+b_k]$ such that the $i$-th block contains $b_i$ elements. If $b_1=\cdots= b_k=s$, we write $\mathcal{OP}_{[s^k]}$ for $\mathcal{OP}_{[b_1,b_2,\ldots,b_k]}$.
 Let $op_{n,k}=|\mathcal{OP}_{n,k}|$, $op_n=|\mathcal{OP}_n|$,  $op_{[b_1,b_2,\ldots,b_k]}=|\mathcal{OP}_{[b_1,b_2,\ldots,b_k]}|$ and $op_{[s^k]}=|\mathcal{OP}_{[s^k]}|$.

 Given an ordered  partition $\pi=B_1/B_2\cdots/B_k\in \mathcal{OP}_{n,k}$ and a
permutation $\sigma=\sigma_1\sigma_2,\ldots,\sigma_m\in S_m$, we say that $\pi$ contains a pattern $\sigma$ if there exist blocks
$B_{i_1},B_{i_2},\ldots, B_{i_m}$ and elements $b_1\in B_{i_1}, b_2\in B_{i_2}, \ldots,
b_m\in B_{i_m}$ with  $1\leq i_1<i_2<\cdots <i_m\leq k$ such that $b_1b_2\cdots b_m$ is order-isomorphic to $\sigma$. Otherwise, we say that $\pi$ avoids a pattern $\sigma$.
For example, the ordered partition $41/53/2\in \mathcal{OP}_{5,3}$ is $123$-avoiding, while it contains a pattern $132$.
Similarly, let $\mathcal{OP}_{n,k}(\sigma)$
denote the set of ordered partitions of $\mathcal{OP}_{n,k}$ that are $\sigma$-avoiding.
Let $op_{n,k}(\sigma)=|\mathcal{OP}_{n,k}(\sigma)|$, $op_n(\sigma)=|\mathcal{OP}_n(\sigma)|$, $op_{[b_1,b_2,\ldots,b_k]}(\sigma)=|\mathcal{OP}_{[b_1,b_2,\ldots,b_k]}(\sigma)|$ and $op_{[s^k]}(\sigma)=|\mathcal{OP}_{[s^k]}(\sigma)|$.

 Godbole, et al. \cite{Pudwell} obtained the following formulas for $op_{n,3}(\sigma)$ and $op_{n,n-1}(\sigma)$ for any $\sigma\in S_3$.

\begin{thm} \label{thm:GOdbole}
For $n\geq 1$, $1\leq k\leq n$, and for any permutation $\sigma$ of length $3$, we have
\begin{align}
op_{n,3}(\sigma)&=\left(\frac{n^2}{8}+\frac{3n}{8}-2\right)2^n+3,\notag\\[3pt]
op_{n,n-1}(\sigma)&=\frac{3(n-1)^2}{n(n+1)}{2n-2 \choose n-1}.\label{equ:opnn-1}
\end{align}
\end{thm}

 Godbole, et al. \cite{Pudwell} also showed  that
\begin{align}
op_{n,k}(\sigma)&=op_{n,k}(123),\label{equ:Godbole}\\[3pt]
op_{[b_1,b_2,\ldots,b_k]}(\sigma)&=op_{[b_1,b_2,\ldots,b_k]}(123)\label{equ:b1b2}
\end{align}
 for any $\sigma\in S_3$. They  raised the question
   concerning the enumeration of $\mathcal{OP}_{n,k}(123)$.
   Using Zeilberger's Maple package $FindRec$ \cite{Findrec}, they
  % posed the following conjecture on $op_{[2^k]}(123)$.
  conjectured that $op_{[2^k]}(123)$ satisfied the following second order linear recurrence relation.

\begin{conj}\label{conj:22}
For $k\geq0$, we have
\begin{align}\label{equ:op2k}
op_{[2^{k+2}]}(123)&=
\frac{329k^3+1215k^2+1426k+528}{2(k+2)(2k+5)(7k+5)}
op_{[2^{k+1}]}(123)\notag\\[6pt]
&\quad \quad +\frac{3(k+1)(2k+1)(7k+12)}{(k+2)(2k+5)(7k+5)}
op_{[2^k]}(123).
\end{align}
\end{conj}

In this paper, we shall give an answer to the above question by providing
 a bivariate generating function for $op_{n,k}(123)$ and we shall
 confirm the conjectured recurrence relation by deriving the
 generating function of $op_{[2^k]}(123)$.

\section{The generating function of $op_{n,k}(123)$}\label{sec:2}

In this section, we obtain the bivariate generating function
of $op_{n,k}(123)$ in answer to the question of Godbole, et al. \cite{Pudwell}.
Let $F(x,y)$ be the generating function of $op_{n,k}(123)$, that is,
\begin{equation}\label{F(x,y)}
F(x,y)=\sum_{n\geq0}\sum_{k\geq0}op_{n,k}(123)x^{n}y^{k}.
\end{equation}

We shall show that $F(x,y)$ can be expressed
in terms of the bivariate generating function $E(x,y)$ of $123$-avoiding
permutations of $[n]$ with respect to the number of descents. More precisely,
  for a permutation $\sigma=\sigma_1\sigma_2\cdots\sigma_n\in S_n$,  the descent set of $\sigma$ is  defined by
  \[ D(\sigma)=\{i:\sigma_i>\sigma_{i+1}\}\]
  and the number of descents of $\sigma$ is denoted by
    $des(\sigma)=|D(\sigma)|$.
Barnabei,  Bonetti and  Silimbani \cite{barnabei} defined the generating function
\begin{equation}\label{E(x,y)_1}
E(x,y)=\sum_{n\geq0}\sum_{\sigma\in S_{n}(123)}x^{n}y^{des(\sigma)}=\sum_{n\geq0}\sum_{d\geq0}e_{n,d}x^{n}y^{d},
\end{equation}
where $e_{n,d}=|\{\sigma~|~\sigma\in S_{n}(123),des(\sigma)=d\}|$. Furthermore, they found the following  formula for $E(x,y)$:
\begin{equation}
E(x,y)=\frac{-1+2xy+2x^{2}y-2xy^{2}-4x^{2}y^{2}+2x^{2}y^{3}
+\sqrt{1-4xy-4x^{2}y+4x^{2}y^{2}}}{2xy^{2}(xy-1-x)}\label{E(x,y)}.
\end{equation}

The main result of this section is stated as follows.

\begin{thm}\label{thm:F(x,y)}
We have
\begin{equation*}
F(x,y)=E(xy,1+y^{-1}),
\end{equation*}
which implies that
  \begin{equation*}
  F(x,y)=\frac{-y-2xy-2x+2x^{2}y+2x^{2}+
     y\sqrt{1-4xy-4x+4x^{2}y+4x^{2}}}{2x(y+1)^{2}(x-1)}.
\end{equation*}
\end{thm}

To prove the above theorem, we establish a
connection between
$op_{n,k}(123)$ and $e_{n,d}$.

\begin{thm}\label{thm1}
For $n\geq 1$ and $1\leq k\leq n$, we have
\begin{equation}\label{op_{n,k}(123)}
op_{n,k}(123)=\sum_{d=n-k}^{n-1}{d \choose n-k}e_{n,d}.
\end{equation}
\end{thm}

\pf Define a map $\varphi\colon\mathcal{OP}_{n,k}(123)\rightarrow S_{n}(123)$
as a canonical representation of an ordered partition.
Given an ordered partition $\pi=B_{1}/B_{2}/\cdots/B_{k}\in\mathcal{OP}_{n,k}(123)$. If we list the elements of each block in decreasing order and ignore the symbol `/'
 between two adjacent blocks, we get   a permutation $\varphi(\pi) =\sigma=\sigma_1\sigma_2\cdots\sigma_n\in S_n$.
We wish to show that $\varphi$ is well-defined, that is, $\sigma=\varphi(\pi)$ is a $123$-avoiding permutation of $S_n$.
Assume to the contrary that $\sigma$ contains a $123$-pattern, that
is, there exist $i < j < l$ such that
 $\sigma_i\sigma_j\sigma_l$ is a $123$-pattern in $\sigma$. By the construction of $\sigma$, we see that the elements $\sigma_i,\sigma_j$ and $\sigma_l$ are in different blocks in $\pi$. This implies that $\sigma_i\sigma_j\sigma_l$ is a $123$-pattern of $\pi$,  a contradiction. Thus   $\sigma\in S_{n}(123)$.
Moreover, according to the construction of $\sigma$, we see that
\begin{equation}\label{ineq:des}
des(\sigma)\geq\sum_{s=1}^{k}(|B_{s}|-1)=n-k.
\end{equation}

Conversely, given a permutation $\sigma=\sigma_{1}\sigma_{2}\cdots\sigma_{n}$ in $S_n(123)$ with $d$ descents, we aim to
count the preimages $\pi$ in $\mathcal{OP}_{n,k}(123)$ such that
$\varphi(\pi)=\sigma$.   If $d< n-k$, by inequality (\ref{ineq:des}),
 it is impossible for any $\pi$ in $\mathcal{OP}_{n,k}(123)$
 to be a preimage of $\sigma$.
 So we may assume that $d\geq n-k$. Let $\pi'=\sigma_{1}/\sigma_{2}/\cdots/\sigma_{n}$. Clearly, $\varphi(\pi')=\sigma$.
%To determine all the preimages $\pi$ of $\sigma$,
If $i\in D(\sigma)$, we may combine $\sigma_i$  and $\sigma_{i+1}$ of $\pi'$ into a block to form a new ordered partition $\pi''$. It is easily
 verified that $\varphi(\pi'')=\sigma$ and $b(\pi'')=n-1$. Moreover,
 we may iterate this process if $des(\pi'')>0$.
Note that at each step we get a preimage of $\sigma$ with one less blocks.
To obtain the preimages $\pi$ with $k$ blocks, we need to repeat this
process $n-k$ times. Observe that the resulting ordered partition depends only on  the positions we choose in $D(\sigma)$. Hence we conclude that
 there are $d\choose n-k$ ordered partitions $\pi$ in $\mathcal{OP}_{n,k}(123)$ such that $\varphi(\pi)=\sigma$.
 This completes the proof.
\qed

Now we are ready to prove Theorem \ref{thm:F(x,y)}.

\noindent
{\it Proof of Theorem \ref{thm:F(x,y)}.}
By Theorem \ref{thm1}, we have
\begin{align*}
\sum_{k\geq0}^{n}op_{n,k}(123)x^{n}y^{k} &=
\sum_{k\geq0}^{n}\sum_{d=n-k}^{n-1}{d \choose n-k}e_{n,d}x^{n}y^{k}\notag\\[3pt]
&=\sum_{d\geq0}^{n-1}\sum_{k=n-d}^{n}{d \choose n-k}e_{n,d}x^{n}y^{k}\notag\\[3pt]
&= \sum_{d\geq0}^{n-1}\sum_{j=0}^{d}{d \choose j}e_{n,d}x^{n}y^{n-j}\notag\\[3pt]
&=\sum_{d\geq0}^{n-1}e_{n,d}(xy)^{n}(1+y^{-1})^{d}.\\[3pt]
\end{align*}
Summing over $n$, we obtain that $F(x,y)=E(xy, (1+y^{-1}))$.
This completes the proof.
\qed

 Setting $y=1$ in the generating function $F(x,y)$, we are led to the
 generating function of $op_n(123)$.

\begin{cor}
Let $H(x)$ be the generating function of $op_{n}(123)$, that is
\begin{align*}
    H(x)=\sum_{n\geq 0}op_n(123)x^n.
\end{align*}
Then we have
\[H(x)=\frac{1}{2}+\frac{1}{1+\sqrt{1-8x+8x^{2}}}.\]
\end{cor}

The connection between  $op_{n,k}(123)$ and $e_{n,d}$ can be
used to derive the following generating function of $op_{n,n-1}(123)$.

\begin{cor}
Let $G(x)$ be the generating function of  $op_{n,n-1}(123)$,
that is,
\[ G(x)=\sum_{n\geq 1}op_{n,n-1}(123)x^{n}.\]
Then we have
\begin{equation}\label{gf:opnn-1}
G(x)={\frac {2\,{x}^{2}-7\,x +2 +3\,x\sqrt {1-4\,x}-2\,\sqrt {1-4\,x}}{2x
\sqrt {1-4\,x}}}.
\end{equation}

\end{cor}

\pf
By Theorem \ref{thm1}, we have
\begin{equation}\label{op_{n,n-1}}
op_{n,n-1}(123)=\sum_{d=1}^{n-1}d e_{n,d}.
\end{equation}
It follows that
\begin{align*}
G(x)&= \sum_{n\geq1}\sum_{d=1}^{n-1}de_{n,d}x^{n}\\[3pt]
&=\frac{\partial E(x,y)}{\partial y}\Big{|}_{y=1}.
\end{align*}

By  expression (\ref{E(x,y)})  for $E(x,y)$,  we obtain (\ref{gf:opnn-1}). This completes the proof. \qed

The formula (\ref{equ:opnn-1}) for $op_{n,n-1}$ can be deduced from  (\ref{gf:opnn-1}).

\section{The generating function of $op_{[2^k]}(123)$}\label{sec:3}

In this section, we derive the generating function of
$op_{[2^k]}(123)$ and confirm  Conjecture \ref{conj:22} on the recurrence relation of
$op_{[2^k]}(123)$.

\begin{thm}\label{thm:Rq}
Let $Q(x)$ be the generating function of $op_{[2^k]}(123)$, that is,
\begin{align*}
    Q(x)=\sum_{k\geq 0}op_{[2^k]}(123)x^{2k}.
\end{align*}
Then we have
\begin{align}\label{gf:Qq}
Q(x)=\sqrt{\frac{2}{1+2\,x^2+\sqrt{1-12\,x^2}}}.
\end{align}
\end{thm}
%\begin{align}\label{gf:Qq}
%Q(x)=\sqrt {-{\frac { \left( 8+2\,{x}^{2} \right)  \left( -2\,{x}^{2}-
%1+\sqrt {-12\,{x}^{2}+1} \right) }{4{x}^{2}}}} \left( 4+{x}^{2}
% \right) ^{-1}.
%\end{align}

%\begin{align}
% Q(x)=\sqrt{\frac{1+2x^2-\sqrt{1-12x^2}}{2x^4+8x^2}}.
%\end{align}

Let $Q^{'}(x)$, $Q^{''}(x)$ and $Q^{'''}(x)$ denote the
first derivative, second derivative and third derivative of $Q(x)$, respectively.
Then it is easily verified that the expression (\ref{gf:Qq})
satisfies the following differential equation.

%which is equivalent to
%the recurrence relation (\ref{equ:op2k}) for $op_{[2^k]}(123)$.

\begin{thm} We have
\begin{align}\label{equ:diffQ}
&\left(\frac{21}{2}x^7+\frac{329}{8}x^5-\frac{7}{2}x^3\right)Q^{'''}(x)
+\left(99x^6+\frac{1443}{8}x^4-5x^2\right)Q^{''}(x)\notag\\[3pt]
&\qquad +\left(207x^5+\frac{717}{8}x^3+11x\right)Q^{'}(x)
+(72x^4-12x^2)Q(x)=0.
\end{align}
\end{thm}

Equating coefficients of $x^{2n+4}$ in   (\ref{equ:diffQ}), we obtain  recurrence relation (\ref{equ:op2k}) for $op_{[2^k]}(123)$.

To prove Theorem \ref{thm:Rq},  we construct a bijection
between ordered partitions and permutations on multisets.
Given an ordered partition $\pi=B_1/B_2/\cdots /B_k\in \mathcal{OP}_{n,k}$,
its canonical sequence, denoted $\psi(\pi)$, is defined to be a sequence $\rho=\rho_1\rho_2\cdots \rho_n$ with $\rho_i=j$ if $i\in B_j$.  Let $\mathcal{W}_{[1^{b_1}2^{b_2}\cdots k^{b_k}]}$ denote the set of
permutations on a multiset $\{1^{b_1},2^{b_2},\ldots ,k^{b_k}\}$,
where $i^r$ means $r$ occurrences of $i$.
It is easily verified that $\psi$ is a bijection between
$\mathcal{OP}_{[b_1,b_2,\ldots,b_k]}$ and $\mathcal{W}_{[1^{b_1}2^{b_2}\cdots k^{b_k}]}$.

Given a permutation $\sigma\in S_m$, if we consider it as an ordered partition
with each block containing only one element, then we can define its canonical sequence to be the canonical sequence of the corresponding ordered partition.
The canonical sequence of $\sigma$ turns out to be the inverse of
$\sigma$, denoted by $\sigma^{-1}$. For example, the canonical sequence of
$43512$ is $45213$.

By the definition of pattern avoiding ordered partitions, we see that
an ordered partition $\pi$ contains a pattern $\sigma$ if and
only if its canonical sequence $\psi(\pi)$ contains a pattern $\sigma^{-1}$.
This implies that  $\psi$ is a bijection between
$\mathcal{OP}_{[b_1,b_2,\ldots,b_k]}(\sigma)$ and $\mathcal{W}_{[1^{b_1}2^{b_2}\cdots k^{b_k}]}(\sigma^{-1})$, where $\mathcal{W}_{[1^{b_1}2^{b_2}\cdots k^{b_k}]}(\tau)$ is the set of
$\tau$-avoiding permutations in  $\mathcal{W}_{[1^{b_1}2^{b_2}\cdots k^{b_k}]}$.
Hence we have
\begin{align}\label{equ:word=op}
    op_{[b_1,b_2,\ldots,b_k]}(\sigma)=|\mathcal{W}_{[1^{b_1}2^{b_2}\cdots k^{b_k}]}(\sigma^{-1})|.
\end{align}

In order to establish recurrence relation for $op_{[2^k]}(123)$, we need
to use $op_{[2^k,1]}(123)$ and $op_{[2^k,1,1]}(123)$.
Combining (\ref{equ:word=op}) and  (\ref{equ:b1b2}), we have
\begin{align*}
    op_{[2^n]}(123)=&|\mathcal{W}_{[1^22^2\cdots n^2]}(132)|,\\[3pt]
    op_{[2^n,1]}(123)=&|\mathcal{W}_{[1^22^2\cdots n^2 (n+1)]}(132)|,\\[3pt]
    op_{[2^{n},1,1]}(123)=&|\mathcal{W}_{[1^22^2\cdots n^2 (n+1) (n+2)]}(132)|.
\end{align*}
Let
\begin{align*}
    u_{2n}&=|\mathcal{W}_{[1^22^2\cdots n^2]}(132)|,\\[3pt]
    u_{2n+1}&=|\mathcal{W}_{[1^22^2\cdots n^2 (n+1)]}(132)|,\\[3pt]
    v_{2n}&=|\mathcal{W}_{[1^22^2\cdots (n-1)^2 n(n+1)]}(132)|,
\end{align*}
we set $u_0=v_0=1$ and set $u_n=v_n=0$ for $n< 0$.

We proceed to derive recurrence relations for $u_{2n}, u_{2n+1}$ and $v_{2n}$
that can be used to obtain a system of equations on the generating functions.
In particular, we get the generating function of $u_{2n}$, that is,
the generating function of $op_{[2^n]}(123)$.

Let $U_e(x)$, $U_o(x)$ and $V(x)$ denote the
generating functions of $u_{2n}$, $u_{2n+1}$ and $v_{2n}$, namely
\begin{align*}
    U_e(x)&=\sum_{n\geq 0}u_{2n}x^{2n},\\[3pt]
    U_o(x)&=\sum_{n\geq 0}u_{2n+1}x^{2n+1},\\[3pt]
    V(x)&=\sum_{n\geq 0}v_{2n}x^{2n}.
\end{align*}

We need the following lemma due to Atkinson, Walker and Linton \cite{Atkinson}.

\begin{lem}\label{symm}
Given two permutations $p=p_1p_2\cdots p_n$ and
$q=q_1q_2\cdots q_n$ on the same multiset of $[n]$,  we have
\begin{align*}
|\mathcal{W}_{[1^{p_1}2^{p_2}\cdots n^{p_n}]}(132)|=
|\mathcal{W}_{[1^{q_1}2^{q_2}\cdots n^{q_n}]}(132)|.
\end{align*}
\end{lem}

\begin{thm}\label{thm:b2n+1}
For $n\geq 0$,
we have
\begin{align}
u_{2n+1}&=\sum_{i+j=2n}u_iu_j,\label{b2n+1}\\[3pt]
U_o(x)&=x\left(U_o^2(x)+U_e^2(x)\right). \label{gf:b_o}
\end{align}
\end{thm}

\pf Assume that $\pi\in \mathcal{W}_{[1^22^2\cdots n^2(n+1)]}(132)$.
Write $\pi$ in the form $\sigma (n+1) \tau$.
Since $\pi$ is $132$-avoiding, both $\sigma$ and $\tau$ are $132$-avoiding. Moreover, for any $r$ in $\sigma$ and
any $s$ in $\tau$, we have $r\geq s$.
Let $k$ be the maximum number in $\tau$.
Then $\tau$ contains all the
numbers in the multiset $\{1^2,2^2,\ldots ,n^2,(n+1)\}$ that are smaller than $k$, that is, $\tau$ contains
a multiset $\{1^2,2^2,\ldots, (k-1)^2\}$.

There are two cases.
If $|\tau|$ is even, then $\tau$  contains two occurrences of $k$.
Hence $\tau$ is in $\mathcal{W}_{[1^22^2\cdots k^2]}(132)$, which is
counted by $u_{2k}$, and $\sigma$ is in $\mathcal{W}_{[(k+1)^2(k+2)^2\cdots n^2]}(132)$. It is easily seen that
$|\mathcal{W}_{[(k+1)^2(k+2)^2\cdots n^2]}(132)|=|\mathcal{W}_{[1^22^2\cdots (n-k)^2]}(132)|$, which is counted by $u_{2n-2k}$.

 If $|\tau|$ is odd, then we have   $\tau\in \mathcal{W}_{[1^22^2\cdots (k-1)^2 k]}(132)$ and $\sigma\in \mathcal{W}_{[k(k+1)^2(k+2)^2\cdots n^2]}(132)$. In this case,  $\mathcal{W}_{[1^22^2\cdots (k-1)^2 k]}(132)$ is counted by $u_{2k-1}$.
  By Lemma \ref{symm}, we see that $|\mathcal{W}_{[k(k+1)^2\cdots n^2]}(132)|=|\mathcal{W}_{[k^2(k+1)^2\cdots (n-1)^2n]}(132)|$,  which is counted by $u_{2n+1-2k}$.
 Combining the above two cases, we obtain  (\ref{b2n+1}).

Using (\ref{b2n+1}), we have
\begin{align*}
    U_o(x)&=\sum_{n\geq0}u_{2n+1}x^{2n+1}\\[3pt]
          &=x\sum_{n\geq0}\sum_{i+j=2n}u_iu_jx^{2n}\\[3pt]
          &=x\sum_{n\geq0}\sum_{2i+2j=2n}u_{2i}u_{2j}x^{2n}
           +x\sum_{n\geq0}\sum_{2i+1+2j+1=2n}u_{2i+1}u_{2j+1}x^{2n}\\[3pt]
          &=x\left(U_o^2(x)+U_e^2(x)\right).
\end{align*}
This completes the proof.\qed

\begin{thm}
For $n\geq 0$,
we have
\begin{align}
v_{2n}&=u_{2n}+u_{2n-1},\label{eq:c_2n}\\[3pt]
V(x)&=U_e(x)+xU_o(x). \label{gf:Cx}
\end{align}
\end{thm}

\pf  Assume that $\pi=\pi_1\pi_2\cdots \pi_{2n} \in \mathcal{W}_{[1^22^2\cdots (n-1)^2n(n+1)]}(132)$.
There are two cases. If $n+1$ precedes $n$ in $\pi$, then
we have $\pi_1=n+1$. Otherwise, $\pi_1(n+1)n$ forms a $132$-pattern of $\pi$, a contradiction. Clearly, in this case  $\pi \in \mathcal{W}_{[1^22^2\cdots (n-1)^2n(n+1)]}(132)$ if and only if $\pi_2\pi_3\cdots \pi_{2n} \in \mathcal{W}_{[1^22^2\cdots (n-1)^2n]}(132)$.  Notice that   $\mathcal{W}_{[1^22^2\cdots (n-1)^2n]}(132)$ is counted by $u_{2n-1}$.

 If $n$ precedes $n+1$ in $\pi$,
     then there does not exist
       any $132$-pattern of $\pi$ that contains both $n$ and $n+1$.
 In this case, we may treat $n+1$ as $n$.
 Such permutations form the set
 $\mathcal{W}_{[1^22^2\cdots (n-1)^2n^2]}(132)$, which is counted by $u_{2n}$.
Combining the above two cases, we obtain (\ref{eq:c_2n}).

The generating function relation (\ref{gf:Cx}) immediately follows
from (\ref{eq:c_2n}). This completes the proof.\qed

\begin{thm}
For $n\geq 1$, we have
\begin{align}
u_{2n}&=2\sum_{2i+j=2n-1}u_{2i}u_j+\sum_{2i+1+j=2n-2}u_{2i+1}u_j-u_{2n-1},
\label{equ:b_2n}\\[9pt]
U_e(x)&=1+2xU_e(x)U_o(x)-x^2U_e^2(x).\label{gf:be}
\end{align}
\end{thm}

\pf  Assume that $\pi\in W_{[1^22^2\cdots n^2]}(132)$.
Write $\pi$ in the form $\sigma n \tau$ such that $\sigma$ contains the other $n$.
Since $\pi$ is $132$-avoiding, both $\sigma$ and $\tau$ are $132$-avoiding. Moreover, for any $r$ in $\sigma$ and
any $s$ in $\tau$, we have $r\geq s$.

Let $k$ be the maximum number in $\tau$.
There are two cases. If $|\tau|$ is even, using the same argument as in
Theorem \ref{thm:b2n+1}, we have $\tau\in \mathcal{W}_{[1^22^2\cdots (k-1)^2k^2]}(132)$ and $\sigma\in \mathcal{W}_{[(k+1)^2\cdots (n-1)^2n]}(132)$.
In this case, $\mathcal{W}_{[1^22^2\cdots (k-1)^2k^2]}(132)$ is counted by
$u_{2k}$ and $\mathcal{W}_{[(k+1)^2\cdots (n-1)^2n]}(132)$ is counted by $u_{2n-1-2k}$.

 If $|\tau|$ is odd,  we have $\tau$  is in $\mathcal{W}_{[1^22^2\cdots (k-1)^2k]}(132)$,
 which is counted by $u_{2k-1}$, and $\sigma$ is in $\mathcal{W}_{[k(k+1)^2\cdots (n-1)^2n]}(132)$. By Lemma \ref{symm}, we see that $|\mathcal{W}_{[k(k+1)^2\cdots (n-1)^2n]}(132)|=
 |\mathcal{W}_{[k^2\cdots(n-2)^2 (n-1)n]}(132)|$,
which is counted by  $v_{2n-2k}$. Notice that $\sigma$ is not empty,
we have $2n-2k>0$.

Combining the above two cases, we get
\begin{align*}
    u_{2n}=\sum_{2i+j=2n-1}u_{2i}u_j+\sum_{2i+1+j=2n-1}u_{2i+1}v_j-u_{2n-1}.
\end{align*}
In view of relation (\ref{eq:c_2n}),  we obtain
\begin{align*}
    u_{2n}&=\sum_{2i+j=2n-1}u_{2i}u_j+\sum_{2i+1+j=2n-1}u_{2i+1}u_j
          +\sum_{2i+1+j=2n-1}u_{2i+1}u_{j-1}-u_{2n-1} \\[6pt]
          &=2\sum_{2i+j=2n-1}u_{2i}u_j+\sum_{2i+1+j=2n-2}u_{2i+1}u_{j}-u_{2n-1}.
\end{align*}

It remains to prove relation (\ref{gf:be}). Using (\ref{equ:b_2n}), we have
\begin{align*}
U_e(x)&=1+\sum_{n\geq 1}u_{2n}x^{2n}\\[3pt]
        &=1+\sum_{n\geq 1}\left( 2\sum_{2i+j=2n-1}u_{2i}u_j+\sum_{2i+1+j=2n-2}
           u_{2i+1}u_{j}-u_{2n-1}\right)x^{2n}\\[9pt]
          &=1+2\sum_{n\geq 1}\sum_{2i+j=2n-1}u_{2i}u_jx^{2n}
          +\sum_{n\geq 1}\sum_{2i+1+j=2n-2}u_{2i+1}u_{j}x^{2n}-
           \sum_{n\geq 1}u_{2n-1}x^{2n}\\[9pt]
          &=1+2xU_e(x)U_o(x)+x^2U_o^2(x)-xU_o(x).
\end{align*}
By equation (\ref{gf:b_o}),  we obtain
\begin{align*}
    U_e(x)&=1+2xU_e(x)U_o(x)+x^2U_o^2(x)-x^2\left(U_o^2(x)+U_e^2(x)\right)\\[6pt]
          &=1+2xU_e(x)U_o(x)-x^2U_e^2(x).
\end{align*}
 This completes the proof.\qed

We are now ready to complete the proof of Theorem \ref{thm:Rq}.

\noindent
{\it Proof of Theorem \ref{thm:Rq}.} Note that $Q(x)=U_e(x)$.
By (\ref{gf:be}), we get
\begin{equation}\label{equ:Bo=Be}
U_o(x)=\frac{x^2U_e^2(x)+U_e(x)-1}{2xU_e(x)}.
\end{equation}
Plugging (\ref{equ:Bo=Be}) into (\ref{gf:b_o}) yields the following equation
\begin{equation}
(x^4+4x^2)U_e^4(x)-(2x^2+1)U_e^2(x)+1=0.
\end{equation}
Given the initial values of $u_{2n}$, we find the  solution of $U_e(x)$
as given by (\ref{gf:Qq}).
 This  completes the proof. \qed

Using (\ref{thm:Rq}), (\ref{gf:Cx}) and (\ref{equ:Bo=Be}),  it can be checked that
\begin{align*}
 U_o(x)&={1\over 2{x}}- {\frac {1+\sqrt {1-12\,{x}^{2}}}{4\,x}}U_e(x),\\[3pt]
V(x)&=\frac{1}{2}+\frac{3-\sqrt {1-12\,{x}^{2}}}{4}U_e(x).
\end{align*}

\vspace{0.5cm}
 \noindent{\bf Acknowledgments.}  This work was supported by  the 973
Project, the PCSIRT Project of the Ministry of Education,  and the National Science
Foundation of China.


\begin{thebibliography}{99}

\bibitem{Atkinson}
M.D. Atkinson, S.A. Linton and L. Walker, Priority queues and multisets, Electron. J. Comb. 2 (1995) \#R24.

%\bibitem{Albert}
%M.H. Albert, R.E.L. Aldred, M.D. Atkinson, C. Handley and D. Holton, Permutations of a multiset avoiding permutations of length 3, European
%J. Combin. 22 (8) (2001) 1021--1031.

\bibitem{barnabei}
M. Barnabei, F. Bonetti and M. Silimbani, The descent statistic on 123-
avoiding permutations, S\'{e}m. Lothar. Combin. 63 (2010) Art. B63a.
%8 pp (electronic).

\bibitem{Pudwell}
A. Godbole, A. Goyt, J. Herdan and L. Pudwell, Pattern avoidance in ordered set partitions, arXiv:math.CO/1212.2530.

\bibitem{A.Goyt}
A. Goyt, Avoidance of partitions of a three-element set, Adv. Appl. Math. 41 (1)  (2008) 95--114.

\bibitem{Mansour}
V. Jel\'{\i}nek and T. Mansour, On pattern-avoiding partitions,
Electron. J. Combin. 15 (2008) \#R39.

\bibitem{M.Klazar}
M. Klazar, On abab-free and abba-free set partitions, European J. Combin. 17 (1) (1996) 53--68.

\bibitem{M.Klazar2}
M. Klazar, Counting pattern-free set partitions I: A generalization of Stirling numbers of the second kind, European J. Combin. 21 (3) (2000) 367--378.

\bibitem{M.Klazar3}
M. Klazar, Counting pattern-free set partitions II: Noncrossing and other hypergraphs, Electron. J. Combin. 7 (2000) \#R34.



\bibitem{Knuth}
D.E. Knuth, The Art of Computer Programming, Sorting and Searching, Vol. 3,
Addison-Wesley, Reading 1973.

\bibitem{B.Sagan}
B.E. Sagan, Pattern avoidance in set partitions, Ars Combin. 94 (2010) 79--96.


\bibitem{opnk(13/2)}
OEIS Foundation Inc., The On-Line Encyclopedia of Integer Sequences, \url{http://oeis.org}.

\bibitem{Findrec}
D. Zeilberger, FindRec: A Maple package that guesses recurrence equations for discrete functions of one variable, \url{http://www.math.rutgers.edu/~zeilberg/EM08/FindRec.txt}, 2008.

\end{thebibliography}
\end{document}